\documentstyle[12pt,amscd]{amsart}
\newtheorem{Theorem}{Theorem}[section] 
\newtheorem{Lemma}[Theorem]{Lemma}
\newtheorem{Corollary}[Theorem]{Corollary}

\newtheorem{Remark}[Theorem]{Remark}

\def\To{\longrightarrow}
\def\reg{\operatorname{reg}}
\def\height{\operatorname{ht}}
\def\im{{\frak m}}
\def\ip{{\frak p}}
\def\ib{{\frak b}}
\def\tor{\operatorname{Tor}}
\def\depth{\operatorname{depth}}
\def\indeg{\operatorname{indeg}}
\def\ass{\operatorname{Ass}}
\def\s{{\sigma}}

\begin{document}

\title{On the regularity of products and intersections \\ of complete intersections}

\author{Marc Chardin}
\address{Institut de Math\'ematiques de Jussieu, CNRS \& Universit\'e Paris VI, Paris, France}
\email{chardin@@math.jussieu.fr}
\author{Nguyen Cong Minh}
\address{Department of Mathematics, University of Education, 136 Xu\^an Thuy, Hanoi, Vietnam}
\email{ngcminhdhsp@@yahoo.com}
\author{Ngo Viet Trung}
\address{Institute of Mathematics, Vi\^en To\`an Hoc, 18 Ho\`ang Qu\^oc Vi\^et, 1037 Hanoi, Vietnam}
\email{nvtrung@@math.ac.vn}
\subjclass{13D02}
\keywords{}
\thanks{The second author is partially supported by the National Basic 
Research Program of Vietnam} 

\begin{abstract}
This paper proves the formulae
\begin{align*}
\reg(IJ) & \le \reg(I) + \reg(J),\\
\reg(I \cap J) & \le \reg(I) + \reg(J)
\end{align*}
for arbitrary monomial complete intersections $I$ and $J$, and provides examples showing that these inequalities do not hold for general  complete intersections. \end{abstract} 

\maketitle

\section{Introduction}

Let $I$ be a homogeneous ideal in a polynomial ring $S$ over a field.
Let
$$0 \To \oplus_jS(-b_{mj}) \To \cdots \To \oplus_jS(-b_{0j}) \To I \To 0$$
be a minimal graded free resolution of $I$. The number
$$\reg(I) := \max_j\{b_{ij}-i|\ i = 0,...,m\}$$
is called the {\it Castelnuovo-Mumford regularity} (or regularity for short) of $I$.
It is of great interest to have good bounds for the regularity [BaM].\par

The regularity of products of ideals was studied first
by Conca and Herzog [CoH]. 
They found some special classes of ideals $I$ and $J$ for which the following formula holds:
$$\reg(IJ) \le \reg(I) + \reg(J)$$
(see also [Si]). In particular, they showed that 
$\reg(I_1\cdots I_d) = d$
for any set of ideals $I_1,...,I_d$ generated by linear forms.
These results led them to raise the question whether the formula
$$\reg(I_1\cdots I_d) \le \reg(I_1) + \cdots + \reg(I_d)$$
holds for any set of complete intersections $I_1,...,I_d$ [CoH, Question 3.6].
Note that this formula does not hold for arbitrary monomial ideals.
For instance, Terai and Sturmfels (see [St]) gave examples of  monomial ideals $I$ such that $\reg(I^2) > 2\reg(I)$.
\par

On the other hand, Sturmfels conjectured that
$\reg(I_1\cap \ldots \cap I_d) \le d$
for any set of ideals $I_1,...,I_d$ generated by linear forms. This conjecture was settled in the affirmative by Derksen and Sidman [DS]. Their proof was inspired by the work of Conca and Herzog. So one might be tempted to ask whether the formula
$$\reg(I_1\cap \cdots \cap I_d) \le \reg(I_1) + \cdots + \reg(I_d)$$
holds for any set of complete intersections $I_1,...,I_d$. \par

The following result show that these question have positive answers in the monomial case and we shall see that there are counter-examples in the general case. 
 
\begin{Theorem}\label{main}
Let $I$ and $J$ be two arbitrary monomial complete intersections. Then
\begin{align*}
\reg(IJ) & \le \reg(I) + \reg(J),\\
\reg(I \cap J) & \le \reg(I) + \reg(J).
\end{align*}
\end{Theorem}

Both formulae follow from a more general bound for the regularity of a larger class of ideals constructed from $I$ and $J$ (Theorem \ref{general 1}).  
The proof is a bit intricate. It is based on a bound for the regularity of a monomial ideal in terms of the degree of the least common multiple of the monomial generators and the height of the given ideal found in [HT].\par

We are not able to extend the first formula to more than two monomial complete intersections. But we find another proof which extends the second formula to any finite set of monomial intersections (Theorem \ref{general 2}).  We would like to mention that the first formula was already proved in the case one of the ideals $I,J$ is generated by two elements by combinatorial methods in [M].\par

In the last section, we give a geometric approach for constructing examples of complete intersection ideals for which the inequalities $\reg(IJ) \le \reg(I) + \reg(J)$ and/or $\reg(I\cap J) \le \reg(I) + \reg(J)$ fails. We show for instance the following:

\begin{Theorem}
Let $Y$ in ${\bf P}^{3}$ be a curve which is defined by at most 4 equations at the generic points of
its irreducible components. Consider 4 elements in $I_{Y}$, $f_{1},f_{2},g_{1},g_{2}$ such that
$I:=(f_{1},f_{2})$ and $J:=(g_{1},g_{2})$ are complete intersection ideals and $I_{Y}$ is the unmixed 
part of $I+J$. Then, if $-\eta :=\min\{ \mu\ \vert\ H^{0}(Y,{\cal O}_{Y}(\mu ))\not= 0\}<0$, one has
$$
\reg (IJ)=\reg (I)+\reg (J)+\eta -1.
$$
\end{Theorem}

A similar construction is explained for $I\cap J$. As a consequence, many families of curves with sections in negative degrees gives rise to counter-examples for the considered inequalities. In the examples we give, $I$ is a monomial ideal and $J$ is either generated by one binomial or by one monomial and one binomial. \smallskip

\noindent{\it Acknowledgment.} The first author is grateful to the Institute of  Mathematics in Hanoi for its hospitality while this work was completed.

\section{Preliminaries}

Let us first introduce some conventions.
For any monomial ideal we can always find a minimal basis consisting of monomials. These monomials will be called the {\it monomial generators} of the given ideal. Moreover, for  a finite set of monomials $A_i = x_1^{a_{i1}}\cdots x_n^{a_{in}}$, we call the monomials $x_1^{\max_i\{a_{i1}\}} \cdots x_n^{\max_i\{a_{in}\}}$  the {\it least common multiple}  of the monomials $A_i$. \par

The key point of our approach is the following bound for the regularity of arbitrary monomial ideals.

\begin{Lemma} \label{HT} {\rm [HT, Lemma 3.1]}
Let $I$ be a monomial ideal. Let $F$ denote the least common multiple of the monomial generators of $I$. Then
$$\reg(I) \le \deg F - \height(I) + 1.$$
\end{Lemma}

This bound is an improvement of the bound $\reg(I) \le \deg F-1$ given by Bruns and Herzog in [BrH, Theorem 3.1(a)].\par

If we apply  Lemma \ref{HT} to the product and the intersection of monomial ideals, we get
\begin{align*}
\reg(I_1\cdots I_d) &\le \sum_{j=1}^d\deg F_j - \height(I_1\cdots I_d)+1,\\
\reg(I_1\cap \cdots \cap I_d) &\le \sum_{j=1}^d\deg F_j - \height(I_1\cap \cdots \cap I_d)+1,
\end{align*}
where $F_j$ denotes the least common multiple of the monomial generators of $I_j$. 
If $I_1,...,I_d$ are complete intersections, then
$\reg(I_j) = \deg F_j - \height(I_j) +1$, whence 
\begin{align*}
\reg(I_1\cdots I_d) &\le \sum_{j=1}^d\reg(I_j)  + \sum_{j=1}^d\height(I_j) - \height(I_1\cdots I_d) - d+1,\\
\reg(I_1\cap \cdots \cap I_d) &\le \sum_{j=1}^d\reg(I_j)  + \sum_{j=1}^d\height(I_j) - \height(I_1\cap\cdots\cap I_d) - d+1.
\end{align*}
These bounds are worse the bounds in the afore mentioned questions. However, the difference is not so big.\par

To get rid of the difference in the case $d=2$ we need the following consequence of Lemma \ref{HT}.

\begin{Corollary} \label{quote}
Let $I$ be a monomial complete intersection and $Q$ an arbitrary monomial ideal (not necessarily a proper ideal of the polynomial ring $S$). Then
$$\reg(I:Q) \le \reg(I).$$
\end{Corollary}

\begin{pf}
Let $F$ denote the product of the monomial generators of $I$. 
Since every monomial generator of $I:Q$ divides a monomial generator of $I$, the least common multiple of the monomial generators of $I:Q$ divides $F$. Applying Lemma \ref{HT} we get
$$\reg(I:Q) \le \deg F - \height(I:Q) + 1 \le \deg F - \height I + 1 = \reg(I).$$
\end{pf}

We will decompose the product and the intersection of two monomial ideals as a sum of smaller ideals and apply the following lemma to estimate the regularity.

\begin{Lemma} \label{sum}
Let $I$ and $J$ be two arbitrary homogeneous ideals. Then
$$\reg (I+J) \le \max\{\reg(I),\reg(J),\reg(I\cap J)-1\}.$$
Moreover, $\reg (I \cap J)= \reg(I + J)+1$ if  $\reg(I+J) > \max\{\reg(I),\reg(J)\}$ or if
$\reg(I\cap J) > \max\{\reg(I),\reg(J)\}+1$.
\end{Lemma}

\begin{pf}
The statements follow from  the exact sequence
$$0 \To I \cap J \To I \oplus J \To I + J \To 0$$ 
and the well-known relationship between regularities of modules of an exact sequence (see e.g. [E, Corollary 20.19]).
\end{pf}

\section{Main results}

We will prove the following general result.

\begin{Theorem}\label{general 1}
Let $I$ and $J$ be two arbitrary monomial complete intersections.
Let $f_1,...,f_r$ be the monomial generators of $I$. Let $Q_1,...,Q_r$ be arbitrary monomial ideals. Then
$$\reg\big(f_1(J:Q_1) + \cdots + f_r(J:Q_r)\big) \le \reg(I) + \reg(J).$$
\end{Theorem}

The formulae of Theorem \ref{main} follow from the above result because
\begin{align*}
IJ  & = f_1J + \cdots + f_rJ = f_1(J:S) + \cdots + f_r(J:S),\\
I \cap J & = (f_1) \cap J + \cdots + (f_r) \cap J = f_1(J:f_1) + \cdots + f_r(J:f_r).
\end{align*}

\begin{pf}
If $r = 1$, we have to prove that
$$\reg \big(f_1(J:Q_1) \big) \le \deg f_1 + \reg(J).$$
It is obvious that
$$\reg\big(f_1(J:Q_1) \big) \le \deg f_1 + \reg(J:Q_1).$$
By Corollary \ref{quote} we have $\reg(J:Q_1) \le \reg(J),$ which implies the assertion. \par

If $r > 1$, using induction we may assume that
\begin{align}
\reg\big(\sum_{i=1}^{r-1}f_i(J:Q_i)\big)  & \le \reg(f_1,...,f_{r-1}),\\
\reg\big(f_r(J:Q_r)\big) & \le \deg f_r + \reg(J).
\end{align}
Since $f_1(J:Q_1),...,f_{r-1}(J:Q_{r-1})$ are monomial ideals, we have 
$$\big(\sum_{i=1}^{r-1}f_i(J:Q_i)\big):f_r = \sum_{i=1}^{r-1}\big(f_i(J:Q_i):f_r\big).$$
Since $f_1,...,f_r$ is a regular sequence, $f_i(J:Q_i):f_r= f_i(J:f_rQ_i)$. 
Therefore,
\begin{align*}
\big(\sum_{i=1}^{r-1}f_i(J:Q_i)\big) \cap f_r(J:Q_r)
& = f_r\big[\big(\sum_{i=1}^{r-1}f_i(J:Q_i):f_r\big) \cap (J:Q_r)\big]\\
& = f_r\big[\big(\sum_{i=1}^{r-1}f_i(J:f_rQ_i)\big) \cap (J:Q_r)\big]\\
& = f_r\big[\sum_{i=1}^{r-1}f_i\big((J:f_rQ_i)\cap (J:f_iQ_r)\big)\big]\\
& = f_r\big[\sum_{i=1}^{r-1}f_i\big(J:(f_rQ_i+f_iQ_r)\big)\big].
\end{align*}
>From this it follows that
$$
\reg\big(\big(\sum_{i=1}^{r-1}f_i(J:Q_i)\big) \cap f_r(J:Q_r)\big) 
 \le \deg f_r + \reg\big(\sum_{i=1}^{r-1}f_i\big(J: (f_rQ_i+f_iQ_r)\big)\big).$$
Using induction we may assume that
$$\reg\big(\sum_{i=1}^{r-1}f_i\big(J:(f_rQ_i+f_iQ_r)\big)\big) \le \reg(f_1,...,f_{r-1}) + \reg(J).$$
Since $\reg(I) = \reg(f_1,...,f_{r-1})+\deg f_r-1$, this implies
\begin{align}
\reg\big(\big(\sum_{i=1}^{r-1}f_i(J:Q_i)\big) \cap f_r(J:Q_r)\big) 
& \le \reg(I) + \reg(J)+1.
\end{align}

Now, we apply Lemma \ref{sum} to the decomposition
$$f_1(J:Q_1)+ \cdots + f_r(J:Q_r) = \big(\sum_{i=1}^{r-1}f_i(J:Q_i)\big)+f_r(J:Q_r).$$
 and obtain
\begin{multline*}
\reg\big(f_1(J:Q_1)+ \cdots + f_r(J:Q_r)\big) \le \\
\max\left\{\reg\big(\sum_{i=1}^{r-1}f_i(J:Q_i)\big),
\reg\big(f_r(J:Q_r)\big),\reg\big(\big(\sum_{i=1}^{r-1}f_i(J:Q_i)\big) \cap f_r(J:Q_r)\big)-1\right\}\\
\le \reg(I) + \reg(J)
\end{multline*}
by using (1), (2), (3).
\end{pf}

\begin{Remark}
{\rm The above proof would work in the case of more than two monomial complete intersections if we 
have a similar result as Lemma \ref{quote}. For instance, if we can prove 
$$\reg(IJ:Q) \le \reg(I) + \reg(J)$$
for two monomial complete intersections $I,J$ and an arbitrary monomial ideal $Q$, then we can give a positive answer to the question of Conca and Herzog in the case $d=3$ for monomial ideals.
We are unable to verify the above formula though computations in concrete cases suggest its validity.}
\end{Remark}

Now we will extend the second formula of Theorem \ref{main} for any set of monomial complete intersections.

\begin{Theorem}\label{general 2}
Let $I_1,\ldots,I_d$ be arbitrary monomial complete intersections. Then
$$\reg(I_1 \cap \cdots \cap I_d) \le \reg(I_1) + \cdots + \reg(I_d).$$
\end{Theorem}

\begin{pf}
We will use induction on the number $n$ of variables and the number 
$$s := \reg(I_1)+ \cdots + \reg(I_d).$$
First, we note that the cases $n = 1$ and $s = 1$ are trivial. \par
Assume that $n \ge 2$ and $r \ge 2$. 
Let $x$ be an arbitrary variable of the polynomial ring $S$. 
It is easy to see that $(I_1,x),...,(I_d,x)$ are monomial complete intersections and
$$(I_1 \cap \cdots \cap I_d, x) = (I_1,x) \cap \cdots \cap (I_d,x).$$
Therefore, using induction on $n$ we may assume that
$$\reg (I_1 \cap \cdots \cap I_d,x) \le \reg(I_1,x) + \cdots + \reg (I_d,x).$$

If $x$ is a non-zerodivisor on $I_1 \cap \cdots \cap I_d$ and if we assume that the intersection is irredundant, then $I_j:x =I_j$ and hence $\reg(I_j,x) = \reg(I_j)$ for all $j = 1,...,d$. In this case,
$$\reg(I_1 \cap \cdots \cap I_d) = \reg(I_1 \cap \cdots \cap I_d, x) \le \reg(I_1)+ \cdots + \reg(I_d).$$

If $x$ is a zerodivisor on $I_1 \cap \cdots \cap I_d$, we involve the ideal
$$(I_1 \cap \cdots \cap I_d):x = (I_1:x) \cap \cdots \cap (I_d:x).$$
If $I_j:x \neq I_j$, either $I_j:x = S$ ($x \in I_j:x$) or $I_j:x$ is a monomial complete intersection generated by the monomials obtained from the generators of $I_j$ by replacing the monomial divisible by $x$ by its quotient by $x$.  In the latter case, we have $\reg(I_j:x) = \reg(I_j)-1$. 
Since there exists at least an ideal $I_j$ with $I_j:x \neq I_j$, the ideal $(I_1 \cap \cdots \cap I_d):x$ is an intersection of monomial complete intersections such that the sum of their regularities is less than $s$. Using induction on $s$ we may assume that
\begin{align*}
\reg((I_1  \cap \cdots \cap I_d):x) &\le \reg(I_1:x) + \cdots + \reg (I_d:x)\\
& \le \reg(I_1)+ \cdots + \reg(I_d)-1.
\end{align*}
Now, from the exact sequence
$$0 \To S/(I_1\cap \cdots \cap I_d):x \overset x \To S/I_1 \cap \cdots \cap I_d \To S/(I_1\cap \cdots \cap I_d,x) \To 0$$ 
we can deduce that
\begin{align*}
\reg(I_1 \cap \cdots \cap I_d) 
& \le \max\left\{\reg((I_1\cap \cdots \cap I_d):x)+1,\reg(I_1\cap \cdots \cap I_d,x)\right\}\\
& \le \reg(I_1) + \cdots + \reg(I_d).
\end{align*}
\end{pf}

\section{Counter-examples}

We will explain a geometric approach, using projective curves, for
constructing families of counter-examples to the inequalities $\reg
(IJ)\leq \reg (I)+\reg (J)$ and $\reg (I\cap J)\leq \reg (I)+\reg
(J)$.  We then give a specific family of such examples, based on
the example [CD, 2.3].

For simplicity, we will work with curves ${\bf P}^{3}$, although this
technique may be easily extended to curves in any projective space.  

Recall that setting, for a finitely generated graded $S$-module $M$,  
$$
a_{i}(M):=\max\{ \mu\ \vert \ H^{i}_{\im}(M)_{\mu}\not= 0\}
$$
if $H^{i}_{\im}(M)\not= 0$ and $a_{i}(M):=-\infty$ else ($\im$ is the
graded maximal ideal of $S$), one has
$$
\reg (M)=\max_{i}\{ a_{i}(M)+i\} .
$$

In the special case were $M=S/I$ is of dimension two (hence defines a
projective scheme of dimension one), one has
$$
\reg (S/I)=\max \{ a_{0}(S/I),a_{1}(S/I)+1,a_{2}(S/I)+2\} ,
$$
and if furthermore $I$ is saturated (in other words is the defining
ideal of the corresponding embedded projective scheme), one has
$a_{0}(S/I):=-\infty$. 

From now on we set $S:=k[x,y,z,t]$ for the homogeneous coordinate ring of ${\bf
  P}^{3}$.\smallskip

\noindent {\bf Step 1.} Construct a curve in ${\bf P}^{3}$ with sections in
negative degrees. \smallskip

This is equivalent to constructing a graded unmixed ideal $I$ 
with $\dim (S/I)=2$, such that
$H^{1}_{\im}(S/I)$ has elements in negative degrees.

One way to construct such a curve is to start from another curve $X$
and two elements in its defining ideal of degrees $d_{1},d_{2}$
that form a complete intersection strictly containing the curve, and
such that the regularity of the ideal of the curve is at least
$d_{1}+d_{2}-1$. 
A good choice is to take a reduced irreducible curve $X$ whose
regularity is at least equal to the sum of the two smallest degree
$d_{1},d_{2}$ of minimal generators of its defining ideal. The
monomial curves offers a good collection of curves of that type.

By liaison, the residual of $X$ in the complete intersection of
degrees $d_{1},d_{2}$ is a non reduced curve $Y$ that  satisfies: 
$$
H^{0}(Y ,{\cal O}_{Y}(\mu ))\not= 0\ \Leftrightarrow
\mu \geq d_{1}+d_{2}-\reg (I_{X})-2
$$
by [CU, 4.2]. In particular $H^{1}_{\im}(S/I_{Y})_{\mu}\not= 0$ for $
 d_{1}+d_{2}-\reg (I_{X})-2\leq \mu <0$ and  
$$
\indeg  (H^{1}_{\im}(S/I_{Y}))=d_{1}+d_{2}-\reg (I_{X})-2. 
$$ 

Also by liaison, $\reg (I_{Y})\leq d_{1}+d_{2}-2$ if $X$ is
reduced. In particular $I_{Y}$ is generated in degrees at most
$d_{1}+d_{2}-2$ in this case.
\smallskip

\noindent {\bf Step 2.} To obtain counter-examples to the inequality $\reg
(I\cap J) \leq \reg (I)+\reg (J)$.\smallskip

Choose three elements in $I_{Y}$
such that they generate an ideal $K :=(f_{1},f_{2},f_{3})$ whose unmixed
part is $I_{Y}$. This is always possible if $Y$ is generically defined by at most
3 equations --in the context of step 1, this is for instance the case if the multiplicities 
of the irreducible components
of $X$ in the complete intersection are at most 3, or if the supports of $X$ and $Y$ are distinct. 
Recall that $I_{Y}$ is generated in degrees at most
$d_{1}+d_{2}-2$ if $X$ is reduced. Therefore if further the supports of
$X$ and $Y$ are distinct (in other words if $Y$ is a geometric link
of the reduced curve $X$ by a complete intersection ideal $\ib$) then
one may choose $f_{1}$ and $f_{2}$  to be the generators of $\ib$ and
$f_{3}\in I_{Y}-\cup_{\ip \in \ass (I_{X})}\ip$ may be chosen of degree at most
$d_{1}+d_{2}-2$.  Then, by [Ch, 0.6],  
$$
\renewcommand{\arraycolsep}{1pt}
\begin{array}{rcl}
a_{0}(S/K)&=&-\indeg ( H^{1}_{\im}(S/I_{Y}))+\s -4,\\
a_{1}(S/K)&=&-\indeg ( I_{Y}/K)+\s -4,\\
a_{2}(S/K)&\leq &-\indeg (K)+\s -5,
\end{array}
$$
where $\s$ is the sum of the degrees of the three forms. 
Since $\indeg ( H^{1}_{\im}(S/I_{Y}))$ is negative, it follows that
$$
\reg (S/K)=a_{0}(S/K)=-\indeg ( H^{1}_{\im}(S/I_{Y}))+\s -4.
$$

Modifying the generators of $K$, if needed, we may assume that
$I=(f_{1},f_{2})$ is a complete intersection ideal and then Lemma \ref{sum}
shows that 
$$
\reg (I\cap (f_{3}))  =\reg (K)+1=-\indeg ( H^{1}_{\im}(S/I_{Y}))+\s -2
 > \s -1=\reg (I)+\reg ((f_{3}))
$$
if (and only if) $\indeg ( H^{1}_{\im}(S/I_{Y}))<-1$.
\smallskip

\noindent {\bf Step 3.} To obtain counter-examples to the inequality $\reg
(IJ) \leq \reg (I)+\reg (J)$. \smallskip

Choose four elements elements
$f_{1},f_{2},g_{1},g_{2}$ in $I_{Y}$ such that:

--- $I:=(f_{1},f_{2})$ and $J:=(g_{1},g_{2})$ are complete
    intersection ideals,

--- $I_{Y}$ is the unmixed part of $I+J$.

Using the example in step 2, one may take the same ideal for $I$, 
$g_{1}:=f_{3}$ and for $g_{2}$ any element in $I_{Y}$ which is prime
to $f_{3}$ (for instance, modifying the generators of $I$, if needed,
one may take $g_{2}:=f_{2}$).

It follows from the isomorphism $\tor_{1}^{S}(S/I,S/J)\simeq (I\cap
J)/IJ$ and the fact that $\depth (S/(I\cap J))>0$ that
$
H^{0}_{\im}(S/IJ)\simeq H^{0}_{\im}(\tor_{1}^{S}(S/I,S/J)).
$
By [Ch, 5.9], $H^{0}_{\im}(\tor_{1}^{S}(S/I,S/J))$ is the graded
$k$-dual of $H^{1}_{\im}(S/I_{Y})$ up to a shift in degrees by $\s'
-4$, where $\s'$ is the sum of the degrees of the 4 forms. We therefore
have
$
H^{0}_{\im}(S/IJ)_{\mu}\simeq H^{1}_{\im}(S/I_{Y})_{\s' -4-\mu}.
$
This implies 
$$
a_{0}(S/IJ)=-\indeg ( H^{1}_{\im}(S/I_{Y}))+\s' -4.
$$

It also follows from the estimates of [Ch, 3.1 (i), (ii) and (iii)] on
$a_{1}(\tor_{1}^{S}(S/I,S/J))$, $a_{2}(\tor_{1}^{S}(S/I,S/J))$ and
$\reg (S/(I+J))$ that 
$$
\reg (IJ)=a_{0}(S/IJ)+1=-\indeg ( H^{1}_{\im}(S/I_{Y}))+\s' -3>\s'
-2=\reg (I)+\reg (J)
$$
if (and only if) $\indeg (H^{1}_{\im}(S/I_{Y}))<-1$. 

Notice that the above considerations shows that any curve $Y$ in ${\bf P}^{3}$ which is
generically defined by at most 4 equations and has sections starting
in degree $-2$ or below gives rise  to counter-examples:

\begin{Theorem} \label{example}
Let $Y$ in ${\bf P}^{3}$ be a curve which is defined by at most 4 equations at the generic points of
its irreducible components. Consider 4 elements in $I_{Y}$, $f_{1},f_{2},g_{1},g_{2}$ such that
$I:=(f_{1},f_{2})$ and $J:=(g_{1},g_{2})$ are complete intersection ideals and $I_{Y}$ is the unmixed 
part of $I+J$. Then, if $-\eta :=\min\{ \mu\ \vert\ H^{0}(Y,{\cal O}_{Y}(\mu ))\not= 0\}<0$, one has
$$
\reg (IJ)=\reg (I)+\reg (J)+\eta -1.
$$
\end{Theorem}

\noindent{\bf A specific class of examples.} 
We consider, as in [CD, 2.3], the monomial curve $Z_{m,n}$
paramaterized on an affine chart by $(1:\theta :\theta^{mn}:\theta^{m(n+1)})$, for
$m,n\geq 2$. 

The binomial $y^{mn}-x^{mn-1}z$
is a minimal generator of the defining ideal of $Z_{m,n}$, hence 
$\reg (I_{Z_{m,n}})\geq mn$. It follows from a theorem of
Bresinski {\it et al.} [BCFH], who determines the regularity of all
curves with a parametrization $(1:t:t^{a}:t^{b})$, that $\reg (I_{Z_{m,n}})=mn$, 
and this may be easily checked in this special case. 

The ideal of this curve contains minimal generators
$x^{m}t-y^{m}z$ and $z^{n+1}-xt^{n}$. 

A component of the scheme defined by the complete intersection ideal
$\ib_{m,n}:=(x^{m}t-y^{m}z,z^{n+1}-xt^{n})$ is the simple line
$x=z=0$, and we take $X:=Z_{m,n}\cup \{ x=z=0 \}$.
On one hand,  $\reg(I_X) \ge mn+1$ because
 $xy^{mn}-x^{mn}z$ is a minimal generator of $I_{X}$. On the other
hand, $\reg (I_{X})\leq \reg (I_{Z_{m,n}})+\reg ((x,z))=mn+1$ because
$\dim (S/I_{Z_{m,n}}+(x,z))\leq 1$. Therefore, $\reg(I_X) = mn+1$.

$Y$ is a geometric link of $X$ with
$$
\indeg (H^{1}_{\im}(S/I_{Y})) = m+n+2-(mn+1)-2=-(m-1)(n-1).  
$$

One has $I_{Y}=(x^{m}t-y^{m}z)+(z,t)^{n}$. To see this, notice that 
$(x^{m}t-y^{m}z)+(z,t)^{n}$ defines a locally complete intersection
scheme supported on the line $z=t=0$, has positive depth and
multiplicity (or degree) at least $n$ --for instance because two
such unmixed ideals differ for two different values of $n$. The
containement    
$$
\ib_{m,n}=I_{X}\cap I_{Y}\subseteq I_{X}\cap
((x^{m}t-y^{m}z)+(z,t)^{n})
$$
and the fact that $\deg (I_{X}\cap I_{Y})=\deg (I_{X})+n$ forces $I_{Y}$ to
coincide with the unmixed ideal $(x^{m}t-y^{m}z)+(z,t)^{n}$. It is
also easy to provide a minimal free $S$-resolution of the ideal
$(x^{m}t-y^{m}z)+(z,t)^{n}$, and show these facts along the same line
as in the proof of [CD, 2.4].\smallskip

For step 2, we take $I:=(t^{n},z^{n})$ and $K:=I+(x^{m}t-y^{m}z)$,
whose saturation is $I_{Y}$, and therefore, 
$$
\reg (I\cap
(x^{m}t-y^{m}z))=(m-1)(n-1)+m+2n-1=(m+1)n
$$
which is bigger than $\reg (I)+\reg
((x^{m}t-y^{m}z))=m+2n$ if and only if $mn>m+n$ ({\it i.e.} iff $(m,n)\not= (2,2)$). 
\smallskip

For step 3, we can take $I:=(t^{n},z^{n})$ and
$J:=(x^{m}t-y^{m}z,t^{n})$, then $\s' =m+3n+1$.
By Theorem \ref{example} we have
$$
\reg(IJ)=(m-1)(n-1)+m+3n-2=mn+2n-1.
$$
Hence,  $\reg (IJ)>\reg (I)+\reg (J)= 
m+3n-1$ if and only if $(m,n)\not= (2,2)$.

\section*{References}

\noindent [BaM] D. Bayer and D. Mumford, What can be computed in algebraic
geometry?  In: D. Eisenbud and L. Robbiano (eds.), Computational Algebraic
Geometry and Commutative Algebra, Proceedings, Cortona 1991, Cambridge
University Press, 1993, 1-48. \par
\noindent [BCFH] H. Bresinsky, F. Curtis, M. Fiorentini, L. T. Hoa, On the structure of local 
cohomology modules for projective monomial curves in ${\bf P}^3$. Nagoya Math. J. 
136 (1994), 81-114.\par
\noindent [BrH] W. Bruns and J. Herzog, On multigraded resolutions, Math.
Proc. Camb. Phil. Soc.  118 (1995), 245-257.\par
\noindent [Ch] M. Chardin, Regularity of ideals and their powers, Pr\'epublication 364, Institut de Mathematiques de Jussieu, 2004. \par
\noindent [CD] M. Chardin, C. D'Cruz, Castelnuovo-Mumford regularity: examples of curves and surfaces, J. Algebra 270 (2003), 347-360. \par
\noindent [CU] M. Chardin, B. Ulrich, Liaison and Castelnuovo-Mumford regularity, Amer. J. Math. 124 (2002),1103-1124.\par
\noindent [CoH] A. Conca and J. Herzog, Castelnuovo-Mumford regularity of products of ideals,
Collect. Math. 54 (2003), 137-152. \par
\noindent [DS] H. Derksen and J. Sidman, A sharp bound for the Castelnuovo-Mumford regularity of subspace arrangements, Adv. Math. 172 (2002), no. 2, 151-157.\par
 \noindent [E] D. Eisenbud, Commutative algebra with a view toward algebraic geometry, Springer, 1994. \par
\noindent [HT] L. T. Hoa and N. V. Trung, On the Castelnuovo-Mumford regularity and the arithmetic degree of monomial ideals,  Math. Zeits. 229 (1998), 519-537. \par
\noindent [M] N. C. Minh, On Castelnuovo-Mumford regularity of products of monomial ideals, Preprint, 2004. \par
\noindent [Si] J. Sidman, On the Castelnuovo-Mumford regularity of products of ideal sheaves, Adv. Geom. 2 (2002), 219-229.\par
\noindent [St] B. Sturmfels, Four counter examples in combinatorial algebraic geometry, J. Algebra 230 (2000), 282-294.\par
\end{document}